\documentclass[a4paper,leqno]{article}
\usepackage{amssymb,amsmath,amsfonts,amsthm,mathrsfs,graphicx,color}
\newtheorem{theorem}{Theorem}[section]

\newtheorem{lemma}[theorem]{Lemma}

\newtheorem{assumption}[theorem]{Assumption}
\newtheorem{remark}[theorem]{Remark}

\begin{document}

\title{\textsf{\textbf{Inverse Problem for the Schr\"{o}dinger Operator\\
in an Unbounded Strip}}} 
\author{Laure Cardoulis\thanks{ Universit\'e de Toulouse 1, UMR 5640, Ceremath/MIP, Place Anatole France, 31000 Toulouse, France, 
laure.cardoulis@univ-tlse1.fr}
\and
Michel Cristofol\thanks{ Universit\'e de Provence, CMI,UMR CNRS 6632, 39, rue Joliot Curie,13453 Marseille Cedex 13
France, Universit\'e Paul C\'ezanne, IUT de Marseille, France, cristo@cmi.univ-mrs.fr}
\and
Patricia Gaitan\thanks{ Universit\'e de Provence, CMI,UMR CNRS 6632, 39, rue Joliot Curie,13453 Marseille Cedex 13, 
France, Universit\'e de la M\'editerran\'ee, IUT d'Aix en Provence, France, gaitan@cmi.univ-mrs.fr}} 
\date{}
\maketitle
\abstract{We consider the operator $H:= i \partial_t + \nabla \cdot (c \nabla)$ 
in an unbounded strip $\Omega$ in $\mathbb{R}^2$, 
where $c(x,y) \in \mathcal{C}^3(\overline{\Omega})$.
We prove an adapted global Carleman estimate and an energy estimate
for this operator. Using these estimates, we give a stability result
for the diffusion coefficient $c(x,y)$.}
%
\section{Introduction}
%
\setcounter{equation}{0}
Let $\Omega =\mathbb{R}  \times (-\frac{d}{2}, \frac{d}{2}) $ be an unbounded strip of $\mathbb{R}^2$
with a fixed width $d$.
We will consider the Schr\"{o}dinger equation
\begin{equation}\label{systq}
\left \{ \begin{array}{ll}
Hq:=i \partial_t q +\nabla \cdot (c(x,y) \nabla q)=0 \;\mbox{ in }\; Q=\Omega \times (0,T),\\
q(x,y,t)=b(x,y,t) \;\mbox{ on }\; \Sigma=\partial\Omega \times (0,T),\\
q(x,y,0)=q_0(x,y) \;\mbox{ on }\; \Omega,
\end{array}\right.
\end{equation}
where $c(x,y) \in \mathcal{C}^3(\overline{\Omega})$ and $c(x,y) \ge c_{min}>0$. Moreover,
we assume that $c$ and all its derivatives up to order three are bounded.
If we assume that $q_0$ belongs to $H^4(\Omega)$ and $b$ is sufficiently regular
(e.g. $b \in H^1(0,T, H^{\frac{9}{2}+\varepsilon}(\partial \Omega)) \cap 
H^2(0,T, H^{\frac{5}{2}+\varepsilon}(\partial \Omega))$ and some additional conditions), then (\ref{systq}) admits a solution
in $H^1(0,T, H^{\frac{3}{2}+\varepsilon}(\Omega))$. We will use this regularity result later.
The aim of this paper is to give a stability and uniqueness result for the coefficient 
$c(x,y)$ using  global Carleman estimates and energy estimates.
We denote by $\nu$ the outward unit normal to $\Omega$ on $\Gamma=\partial \Omega$. 
We denote $\Gamma=\Gamma^+ \cup \Gamma^-$, where
$\Gamma^+=\{(x,y) \in \Gamma ; \; y=\frac{d}{2}\}$ and 
$\Gamma^-=\{(x,y) \in \Gamma ; \; y=-\frac{d}{2}\}$. 
We use the following notations
$\nabla \cdot (c \nabla u)= \partial_x(c \partial_x u)+ \partial_y(c \partial_y u)$,
$\nabla u \cdot \nabla v= \partial_x u \partial_x v + \partial_y u \partial_y v$,
$\partial_{\nu} u=\nabla u \cdot \nu$. 
\\ \noindent
We shall use the following notations $Q=\Omega \times (0,T)$, $\widetilde{Q}=\Omega \times (-T,T)$,
$\Sigma=\Gamma \times (0,T)$, $\widetilde{\Sigma}=\Gamma \times (-T,T)$, 
$\Lambda(R_1):= \{\Phi \in L^{\infty}(\Omega), 0<R_{1} \leq \|\Phi\|_{L^{\infty}(\Omega)} \},$
and
$\Lambda(R_2):= \{\Phi \in L^{\infty}(\Omega), \|\Phi\|_{L^{\infty}(\Omega)} \leq R_2 \},$
where $R_{1}$ and $R_{2}$ are positive constants with $R_1 \leq R_2$.\\ \\ \noindent
Our problem can be stated as follows:\\  \noindent
Is it possible to determine the coefficient $c(x,y)$ from the measurement of
$\partial_{\nu}(\partial_t q)$ on $\Gamma^+$?\\ \\ \noindent
Let $q$ (resp. $\widetilde{q}$) be a solution of (\ref{systq}) associated with
($c$, $b$, $q_0$) (resp. ($\widetilde{c}$, $b$, $q_0$)) satisfying
some regularity properties:
\begin{itemize}
\item $\partial_t \widetilde{q}$, $\nabla (\partial_t \widetilde{q})$ and 
$\Delta(\partial_t \widetilde{q})$ are in $\Lambda(R_2)$,
\item $q_0$ is a real valued function in $\mathcal{C}^3(\Omega)$,
\item $q_0$ and all its derivatives up to order three are in $\Lambda(R_2)$ .
\end{itemize}

\noindent 

Our main result is
$$|c-\widetilde{c}|^2_{H^1(\Omega)} \leq C|\partial_{\nu}(\partial_t q)-
\partial_{\nu}(\partial_t \widetilde{q})|^2_{L^2((0,T)\times \Gamma^+)},$$
where $C$ is a positive constant which depends on $(\Omega, \Gamma, T, R_1, R_2)$
and where the above norms are weighted Sobolev norms. \\
 \noindent
 The major novelty of this paper is to give an $H^1$ stability estimate
 for the diffusion coefficient with only one observation in an unbounded domain.\\
 \noindent
We prove an adapted global Carleman estimate and an energy estimate 
for the operator $H$ with a boundary term on $\Gamma^+$. Such energy estimate
has been proved in \cite{LTZ} for the Schr\"{o}dinger operator in a
bounded domain in order to obtain a controllability result.
Then using these estimates and following the method developed by
Imanuvilov, Isakov and Yamamoto for the Lam\'e system in \cite{IIY}, \cite{IY}, we give
a stability and uniqueness result for the diffusion coefficient $c(x,y)$.
Note that this stability result corresponds to a stability result for three linked coefficients 
($c$, $\partial_y c$ and $\partial_y c$) with only one observation. For independent coefficients, in our knowledge,
there is no stability result with one observation. \\
\noindent
The method of Carleman estimates was introduced in the field of inverse problems in the
works of Bukhgeim and Klibanov (see \cite{B:99}, \cite{BK:81}, \cite{K:84}, \cite{K:92}).
The first stability result for a multidimensional inverse problem (for a hyperbolic equation) was
obtained by Puel and Yamamoto \cite{PY:96} using a modification of the idea of \cite{BK:81}.\\ \noindent
For the non stationnary Schr\"{o}dinger equation, \cite{BP} gives a stability result for the
potential in a bounded domain.
For the stationnary Schr\"{o}dinger equation, we can cite recent results concerning uniqueness
for the potential from partial Cauchy data (see for exemple \cite{KSU} and the references herein).
\\
\noindent 
In unbounded domains Carleman estimate with an internal observation has been proved for 
the heat equation  in \cite{BT}.
 \\ \noindent
A physical background could be the characterization of the diffusion 
coefficient for a strip in geophysics. 
Indeed if we look for time harmonic solutions of (\ref{systq}), the problem can be written,
after some changes of variables as the 
reconstruction of a non local potential $P$ in a strip for the operator $-\Delta +P$.
Few results for inverse problems exist in a two-dimensional strip (see \cite{CG}).
For the layer $\mathbb{R}^n \times [0,h]$ with $n\geq 2,$ several results exist 
for the stationnary inverse problems (see \cite{BGWX}, \cite{CGI}, \cite{CW}, \cite{I}, 
 \cite{IMN}, \cite{W}, ...).\\ \noindent
On the other hand, we can link our problem to the determination 
of the curvature function for a curved quantum guide 
(see \cite{ES}, \cite{CDFK}, \cite{DEK}, ...).\\  
\noindent 
This paper is organized as follows. 
In section $2$, we give an adapted global Carleman estimate for the 
operator $H$.
In section $3$, we prove an energy estimate and we give a stability 
result for the diffusion coefficient $c.$
%
\section{Global Carleman Estimate}
%
\setcounter{equation}{0}
Let $c=c(x,y)$ be a bounded positive function in $\mathcal{C}^3(\overline{\Omega})$ such that
\begin{assumption} 
\label{c}
$c(x,y) \in \Lambda(R_1)$, $c$ and all its derivatives up to order three are in $\Lambda(R_2)$.
\end{assumption}
Let $q=q(x,y,t)$ be a function equals to zero on $\partial \Omega \times (-T,T)$ and solution
of the Schr\"{o}dinger equation 
$$i \partial_t q +\nabla \cdot (c(x,y) \nabla q)=f.$$
We prove here a global Carleman-type estimate for $q$
with a single observation 
acting on a part $\Gamma^+$ of the boundary $\Gamma$ 
in the right-hand side of the estimate.  
Let $\widetilde{\beta}$ be a 
${\cal C}^4(\overline{\Omega})$ positive function such that there exists positive
constant $C_{pc}$ which satisfies 
\begin{assumption}
\label{funct-beta}
\begin{itemize}
\item
$|\nabla \widetilde{\beta}| \in \Lambda(R_1), \;\;
{\partial}_{\nu} {\widetilde{\beta}}\leq 0\;\;\mbox{on}\;\;\Gamma^-$,
\item
$ \widetilde{\beta}$ and all its derivatives up to order four are in $\Lambda(R_2)$.
\item
$2 \Re (D^2 \widetilde{\beta}(\zeta, \bar{\zeta})) 
-c \nabla c \cdot \nabla \widetilde{\beta} |\zeta|^2+2c^2|\nabla \widetilde{\beta} \cdot \zeta|^2
 \geq C_{pc} |\zeta|^2$, for all $\zeta \in \mathbb{C}$
\end{itemize}
\end{assumption}
where 
$$D^2 \widetilde{\beta}=\left( \begin{array}{cc}
c\partial_x(c\partial_x \widetilde{\beta}) & c\partial_x(c\partial_y \widetilde{\beta})\\
c\partial_y(c\partial_x \widetilde{\beta}) & c\partial_y(c\partial_y \widetilde{\beta})
\end{array} \right).$$
Note that the last assertion of Assumption \ref{funct-beta} expresses
the pseudo-convexity condition for the function $\widetilde{\beta}$.
This Assumption imposes restrictive conditions for the choice of the functions
$c(x,y)$ in connection with the function $\widetilde{\beta}$.
Note that there exists functions satisfying such Assumptions; indeed, if we consider
$$c(x,y) \in \left\{ f \in C^1(\Omega); \exists r_0 \mbox{ positive constant},
\left\{ \begin{array}{ll} 
-f\ \partial_y f \partial_y \widetilde{\beta} \geq r_0>0,\\ 
f\ \partial_y f \partial_y \widetilde{\beta} ((\frac{\partial_x f}{\partial_y f})^2+1)
+2f^2(\partial_{yy} \widetilde{\beta}+(\partial_y \widetilde{\beta})^2)\geq r_0>0.
\end{array}\right. \right\}$$
then a function $\widetilde{\beta}(x,y)=\widetilde{\beta}(y)$ is available 
(for example, $c(x,y)=(\frac{1}{1+x^2}+1) e^{-y}$ and $\widetilde{\beta}(x,y)=e^{y}$).\\[2mm]
Similar restrictive conditions have been highlighted for the hyperbolic case
in \cite{KY}, \cite{KT:04} and for the Schr\"{o}dinger operator in \cite{I:98}\\
\noindent
Then, we define $\beta= \widetilde{\beta}+K$ with
$K= m \|\widetilde{\beta}\|_{\infty}$ and $m>1$. For $\lambda> 0$ and $t \in
(-T,T)$, we define the following weight functions
\begin{equation}
  \label{wf}
  \varphi(x,y,t)=\frac{e^{\lambda \beta(x,y)}}{(T+t)(T-t)},
  \quad \quad \eta(x,y,t)=\frac{e^{2\lambda K} -e^{\lambda
  \beta(x,y)}}{(T+t)(T-t)}.
\end{equation}
Let $H$ be the operator defined by
\begin{equation} \label{H}
Hq:=i\partial_t q +\nabla \cdot (c(x,y) \nabla q)  \;\mbox{ in }\; \widetilde{Q}=\Omega \times (-T,T).
\end{equation}
We set $\psi=e^{-s \eta}q$, $M \psi = e^{-s \eta} H(e^{s \eta} \psi)$ for $s>0$
 and we introduce the following operators
\begin{equation} \label{M1}
  M_1\psi  : = i\partial_t \psi+\nabla \cdot (c \nabla\psi)
+s^2 c |\nabla \eta |^2 \psi,
\end{equation}
\begin{equation} \label{M2}
  M_2\psi  : = is \partial_t\eta \psi+2 c s \nabla \eta  \cdot \nabla \psi
+s \nabla \cdot (c \nabla \eta) \psi.
\end{equation}
Then the following result holds.\\
\begin{theorem}
\label{th-Carl} 
Let $H$, $M_1$, $M_2$ be the operators defined respectively by
(\ref{H}), (\ref{M1}), (\ref{M2}). We assume that Assumptions \ref{c} and \ref{funct-beta}
are satisfied. 
Then there exist $\lambda_0> 0$, $s_0>0$ and a positive
constant $C=C(\Omega, \Gamma,T, C_{pc}, R_1, R_2)$ such that, for any $\lambda \ge
\lambda_0$ and any $s \ge s_0 $, the next inequality holds:
\begin{eqnarray}
\label{Carl}
   s^{3} \lambda^{4}\int_{-T}^T   \int_{\Omega} e^{-2s \eta} \varphi^{3} |q|^2\ d x \ dy \ d t  
+s \lambda \int_{-T}^T   \int_{\Omega} e^{-2s \eta} \varphi |\nabla q|^2\ d
x \ dy \ d t
  +  \|M_1(e^{-s \eta}q)\|^2_{L^2(\widetilde{Q})} \\
+ 
\|M_2(e^{-s \eta}q)\|^2_{L^2(\widetilde{Q})}
  \leq C \left[
   s \lambda \int_{-T}^T  
\int_{\Gamma^+} e^{-2s \eta} \varphi  |\partial_{\nu} q|^2\ \partial_{\nu} \beta \ d \sigma \ d t
  +\int_{-T}^T   \int_{\Omega} e^{-2s \eta}\ | H q |^2\ d x \ dy \ d t\right],\nonumber
\end{eqnarray}
for all $q$ satisfying $Hq\in L^2(\Omega \times (-T,T)),$ $q\in
L^2(-T,T;H^1_0(\Omega)),$ $\partial_{\nu} q\in L^2(-T,T;L^2(\Gamma)).$
\end{theorem}
\begin{proof}
If we set $\psi=e^{-s \eta}q$, we calculate $M\psi=e^{-s \eta}H(e^{s \eta} \psi)$ and we obtain:
$$M\psi=M_1 \psi + M_2 \psi $$
with $M_1$ and $M_2$ defined respectively by (\ref{M1}) and (\ref{M2}).
Then 
\begin{eqnarray} \label{eqps}
\int_{-T}^T  \int_{\Omega}|M\psi|^2dx \ dy \ dt& = &
\int_{-T}^T   \int_{\Omega}|M_1\psi|^2dx \ dy \ dt
+\int_{-T}^T   \int_{\Omega}|M_2\psi|^2dx \ dy \ dt\\ \nonumber
& + & 2 \Re (\int_{-T}^T   \int_{\Omega}M_1\psi\  \overline{M_2\psi}\ dx \ dy \ dt),
\end{eqnarray}
where $\overline{z}$ is the conjugate of $z$, $\Re \ (z)$ its real part
and $\Im \ (z)$ its imaginary part.
We have to compute the scalar product in (\ref{eqps})
$$\Re (\int_{-T}^T   \int_{\Omega}M_1\psi \ \overline{M_2\psi}\ dx \ dy \ dt)=
I_{11'}+I_{12'}+I_{13'}+I_{21'}+I_{22'}+I_{23'}+I_{31'}+I_{32'}+I_{33'}.$$
Then, we have
\begin{equation} \label{I11}
I_{11'}=\Re \left( \int_{-T}^T  \int_{\Omega} (i\partial_t \psi)
(-is\partial_t\eta \ \overline{\psi})\  dx \ dy \ dt \right)=
-\frac{s}{2} \int_{-T}^T  \int_{\Omega} \partial_{tt} \eta\  |\psi|^2 \ dx \ dy \ dt.
\end{equation}
\begin{eqnarray*}
I_{12'} & = & \Re\left ( 2i  s \int_{-T}^T  \int_{\Omega}
 c\ \partial_t \psi \nabla \eta \cdot \nabla \overline{\psi} \ dx \ dy \ dt\right)\\
& = & s \Im \left(\int_{-T}^T  \int_{\Omega} 
c\ \nabla \eta \cdot \nabla \psi \ \partial_t \overline{\psi} \ dx \ dy \ dt\right)
- s\Im \left(\int_{-T}^T  \int_{\Omega} 
c\ \nabla \eta \cdot \nabla \overline{\psi} \  \partial_t \psi \ dx \ dy \ dt \right).
\end{eqnarray*}
After an integration by parts with respect to the space variable in the first integral
and to the time variable in the second integral, we obtain
\begin{eqnarray} \label{I12}
I_{12'} =  - s \Im \left( \int_{-T}^T  \int_{\Omega}
\nabla \cdot (c \nabla \eta)\  \psi\ \partial_t \overline{\psi}\  
 dx \ dy \ dt \right)
+   s \Im \left(\int_{-T}^T  \int_{\Omega}
c \psi\ \nabla( \partial_t \eta) \cdot \nabla \overline{\psi}\ 
dx \ dy \ dt\right).
\end{eqnarray}
\begin{eqnarray} \label{I13}
I_{13'} & = & s \Im \left(\int_{-T}^T  \int_{\Omega}
\nabla \cdot(c \nabla \eta)\ \psi \ \partial_t \overline{\psi}\  
dx \ dy \ dt\right).
\end{eqnarray}
Note that $I_{13'}$ vanishes with the first term of $I_{12'}$.
\begin{eqnarray} \label{I21}
I_{21'} & = & \Re \left(-is \int_{-T}^T  \int_{\Omega} \partial_t \eta \ \overline{\psi}\  
\nabla \cdot(c \nabla \psi)  dx \ dy \ dt \right)\\ \nonumber
& = & s \Im \left(\int_{-T}^T  \int_{\Omega}
c \psi\ \nabla (\partial_t \eta) \cdot \nabla \overline{\psi}\ dx \ dy \ dt\right).
\end{eqnarray}
Using integrations by parts, we obtain
\begin{eqnarray}
I_{22'} & = & 2 s  \Re \left(\int_{-T}^T  \int_{\Omega} 
c \nabla \cdot ( c \nabla \psi) \nabla \eta \cdot \nabla \overline{\psi}) dx \ dy \ dt\right) \\ \nonumber
& = & - s \lambda \int_{-T}^T  \int_{\Omega}  \varphi 
| \nabla \psi |^2 ( \nabla \cdot (c^2 \nabla \beta) + \lambda c^2 |\nabla \beta|^2)
\ dx \ dy \ dt \\ \nonumber
& + &  s  \int_{-T}^T  \int_{\partial \Omega} c^2 \partial_{\nu} \eta
\ |\partial_{\nu} \psi |^2 \ d\sigma \ dt
 +  2 s \lambda^2  \int_{-T}^T  \int_{\Omega}
\varphi c^2 | \nabla \beta \cdot \nabla \psi|^2 dx  \  dy \ dt\\ \nonumber
& + & 2 s \lambda  \Re \left(\int_{-T}^T  \int_{ \Omega}
\varphi c  \sum_{i,j=1}^2 \partial_{x_i}(c \partial_{x_j} \beta)\ \partial_{x_i} \psi
\ \partial_{x_j} \overline{\psi}\ dx \ dy  \ dt\right) \\ \nonumber
& = & - s \lambda \int_{-T}^T  \int_{\Omega}  \varphi 
| \nabla \psi |^2 ( \nabla \cdot (c^2 \nabla \beta) + \lambda c^2 |\nabla \beta|^2)
\ dx \ dy \ dt \\ \nonumber
& - &  s \lambda \int_{-T}^T  \int_{\partial \Omega} c^2 \varphi\ \partial_{\nu} \beta
\ |\partial_{\nu} \psi |^2 \ d\sigma \ dt
 +  2 s \lambda^2  \int_{-T}^T  \int_{\Omega}
\varphi c^2 | \nabla \beta \cdot \nabla \psi|^2 dx  \  dy \ dt\\ \nonumber
& + & 2 s \lambda  \Re \left(\int_{-T}^T  \int_{ \Omega}
\varphi D^2 \beta (\nabla \psi, \nabla \overline{\psi})\ dx \ dy  \ dt\right).
\end{eqnarray}
Using integration by parts, we obtain
\begin{eqnarray}
I_{23'} & = &  s  \Re\left(  \int_{-T}^T  \int_{\Omega} 
\nabla \cdot(c \nabla \psi)\ \nabla \cdot ( c \nabla \eta)
\ \overline{\psi}\ dx \ dy \ dt\right)\\ \nonumber
& = &s \lambda \int_{-T}^T  \int_{\Omega} c \varphi 
\ | \nabla \psi |^2 ( \nabla \cdot (c \nabla \beta) + \lambda c 
|\nabla \beta|^2)\ dx \ dy \ dt \\ \nonumber
& - & \frac{s \lambda^2}{2} \int_{-T}^T  \int_{\Omega}
\varphi \nabla \cdot(c \nabla \cdot(c\nabla \beta)\nabla \beta)) |\psi|^2 dx
\ dy \ dt \\ \nonumber
& - & \frac{s \lambda^3}{2}\int_{-T}^T  \int_{\Omega}
\varphi c |\nabla \beta|^2 \nabla \cdot( c \nabla \beta) |\psi|^2 dx \ dy \ dt \\ \nonumber
& - & \frac{s \lambda^4}{2}\int_{-T}^T  \int_{\Omega}
\varphi c^2 |\nabla \beta|^4 |\psi|^2 dx \ dy \ dt
-\frac{s \lambda^3}{2}\int_{-T}^T  \int_{\Omega}
\varphi \nabla \cdot(c^2 |\nabla \beta|^2 \nabla \beta) |\psi|^2 dx \ dy \ dt \\ \nonumber
& - & \frac{s \lambda^2}{2}\int_{-T}^T  \int_{\Omega}
\varphi c \nabla \beta \cdot \nabla  ( \nabla \cdot (c \nabla \beta) + \lambda c |\nabla \beta|^2)
 |\psi|^2 dx \ dy \ dt \\ \nonumber
 & - & \frac{s \lambda}{2}\int_{-T}^T  \int_{\Omega}
\varphi  \nabla \cdot ( c \nabla  ( \nabla \cdot (c \nabla \beta) + \lambda c |\nabla \beta|^2))
 |\psi|^2 dx \ dy \ dt.
\end{eqnarray}
And we obviously have
\begin{equation} \label{I31}
I_{31'}=s^3 \Re \left(  \int_{-T}^T  \int_{\Omega}c\ (-i \partial_t \eta \ \overline{\psi})
|\nabla \eta|^2 \psi \ dx \ dy \ dt\right) =0.
\end{equation}
\begin{eqnarray}
I_{32'}& = & 2 s^3  \Re \left(\int_{-T}^T  \int_{\Omega} 
c^2 |\nabla \eta|^2 \psi\  \nabla \eta \cdot \nabla \overline{\psi}\ dx \ dy \ dt\right)\\ \nonumber 
& =& s^3 \int_{-T}^T  \int_{\Omega} 
c^2 |\nabla \eta|^2 \nabla \eta \cdot \nabla |\psi|^2\ dx \ dy \ dt
\\ \nonumber
 & = & s^3 \lambda^3 \int_{-T}^T  \int_{\Omega}
c \varphi^3 \left( |\nabla \beta|^2 \nabla \cdot (c \nabla \beta) +
 \nabla \beta \cdot \nabla(c|\nabla \beta|^2) \right) |\psi|^2 dx \ dy \ dt \\ \nonumber
 & + & 3 s^3 \lambda^4 \int_{-T}^T  \int_{\Omega}
c^2 \varphi^3 |\nabla \beta|^4 |\psi|^2 dx \ dy \ dt .
\end{eqnarray}
\begin{eqnarray} \label{I33}
I_{33'} & = & s^3 \Re\left(  \int_{-T}^T  \int_{\Omega} 
c\ |\nabla \eta|^2 \psi \nabla \cdot ( c \nabla \eta) \overline{\psi}
\ dx \ dy \ dt\right) \\ \nonumber
  & = &-s^3 \lambda^3 \int_{-T}^T  \int_{\Omega}
c\ \varphi^3 |\nabla \beta|^2 \nabla \cdot(c \nabla \beta) |\psi|^2
dx \ dy \ dt \\ \nonumber
&- & s^3 \lambda^4 \int_{-T}^T  \int_{\Omega}
 c^2 \varphi^3 |\nabla \beta|^4 |\psi|^2
dx \ dy \ dt.  
\end{eqnarray}
So, by (\ref{I11})-(\ref{I33}), we get :
\begin{eqnarray*}
\Re\left( \int_{-T}^T   \int_{\Omega}M_1\psi\ \overline{M_2\psi}\ dx \ dy \ dt\right)&=&
2 s \lambda^2  \int_{-T}^T  \int_{\Omega}
\varphi c^2 | \nabla \beta \cdot \nabla \psi|^2 dx  \  dy \ dt \\
& - & s \lambda \int_{-T}^T  \int_{\Omega}
\varphi\ c \ \nabla c \cdot \nabla \beta |\nabla \psi|^2 dx  \  dy \ dt \\
& + & 2 s \lambda  \Re \left(\int_{-T}^T  \int_{ \Omega}
\varphi D^2 \beta(\nabla \psi, \nabla \overline{\psi})\  dx \ dy  \ dt\right) \\
& + & 2 s^3 \lambda^4 \int_{-T}^T  \int_{\Omega}
c^2 \varphi^3 |\nabla \beta|^4 |\psi|^2 dx \ dy \ dt \\
& - &  s \lambda \int_{-T}^T  \int_{\partial \Omega} c^2 \varphi \partial_{\nu} \beta
|\partial_{\nu} \psi |^2 \ d\sigma \ dt + X,
\end{eqnarray*}
where 
\begin{eqnarray*}
X & = & \frac{-s}{2} \int_{-T}^T  \int_{\Omega} \partial_{tt} \eta\ |\psi|^2\ dx \ dy \ dt
+ 2 s \Im \int_{-T}^T  \int_{\Omega}
c \psi \nabla( \partial_t \eta) \cdot \nabla \overline{\psi}
dx \ dy \ dt \\
& - & \frac{s \lambda^2}{2} \int_{-T}^T  \int_{\Omega}
\varphi \nabla \cdot(c \nabla \cdot(c\nabla \beta) \nabla \beta)) |\psi|^2 dx \ dy \ dt \\
& -& \frac{s \lambda^3}{2}\int_{-T}^T  \int_{\Omega}
\varphi c |\nabla \beta|^2 \nabla \cdot( c \nabla \beta) |\psi|^2 dx \ dy \ dt \\ 
& - & \frac{s \lambda^4}{2}\int_{-T}^T  \int_{\Omega}
\varphi c^2 |\nabla \beta|^4 |\psi|^2 dx \ dy \ dt
-\frac{s \lambda^3}{2}\int_{-T}^T  \int_{\Omega}
\varphi \nabla \cdot(c^2 |\nabla \beta|^2 \nabla \beta) |\psi|^2 dx \ dy \ dt \\
& - & \frac{s \lambda^2}{2}\int_{-T}^T  \int_{\Omega}
\varphi c \nabla \beta \cdot \nabla  ( \nabla \cdot (c \nabla \beta) + \lambda c |\nabla \beta|^2)
 |\psi|^2 dx \ dy \ dt \\ 
 & - & \frac{s \lambda}{2}\int_{-T}^T  \int_{\Omega}
\varphi  \nabla \cdot ( c \nabla  ( \nabla \cdot (c \nabla \beta) + \lambda c |\nabla \beta|^2))
 |\psi|^2 dx \ dy \ dt\\
 & + & s^3 \lambda^3 \int_{-T}^T  \int_{\Omega}
\varphi^3 
c \nabla \beta \cdot \nabla(c|\nabla \beta|^2) |\psi|^2 dx \ dy \ dt.
\end{eqnarray*}
Recall that:\\
$$ 
\|M\psi\|^2_{L^2(\widetilde{Q})} = \|M_1\psi\|^2_{L^2(\widetilde{Q})} +
 \|M_2\psi\|^2_{L^2(\widetilde{Q})} + 2 \Re (M_1\psi,\overline{M_2 \psi}).
$$
Then :
\begin{eqnarray} \label{ine}
\|M_1\psi\|^2_{L^2(\widetilde{Q})}  +  \|M_2\psi\|^2_{L^2(\widetilde{Q})} + 4 s \lambda^2  \int_{-T}^T  \int_{\Omega}
\varphi c^2 | \nabla \beta \cdot \nabla \psi|^2 dx  \  dy \ dt \\ \nonumber
 +  4 s \lambda  \Re \left(\int_{-T}^T  \int_{ \Omega} 
\varphi D^2 \beta(\nabla \psi, \nabla \overline{\psi}) dx \ dy  \ dt\right) 
-2s \lambda \int_{-T}^T  \int_{\Omega}
\varphi\ c \ \nabla c \cdot \nabla \beta |\nabla \psi|^2 dx  \  dy \ dt \\ \nonumber
 +  4 s^3 \lambda^4 \int_{-T}^T  \int_{\Omega}
c^2 \varphi^3 |\nabla \beta|^4 |\psi|^2 dx \ dy \ dt 
-   2s \lambda \int_{-T}^T  \int_{\partial \Omega} c^2 \varphi \partial_{\nu} \beta
|\partial_{\nu} \psi |^2 \ d\sigma \ dt 
 \le  \|M\psi\|^2_{L^2(\widetilde{Q})} + 2 |X|.
\end{eqnarray}
Taking into account
\begin{itemize}
\item
$|\widetilde{\beta}|+|\nabla \widetilde{\beta}|+
|\nabla(\nabla \cdot (c\nabla \widetilde{\beta}))|+
|\nabla \cdot (\nabla(\nabla \cdot (c\nabla \widetilde{\beta})))|
\leq C(\Omega, \Gamma, T, R_2)
 \;\;\mbox{ in } \;\;\Omega,$
\item
$|\partial_{tt} \eta| \leq C(T) \varphi^3,$
\quad
$|\partial_t \varphi| \leq C(T) \varphi^2,$
\quad
$\varphi \leq C(T) \varphi^3,$
\quad
$\varphi^2 \leq C(T) \varphi^3,$
\item
${\displaystyle{ | s \Im ( \int_{-T}^T  \int_{\Omega}
c \psi \nabla( \partial_t \eta) \cdot \nabla \overline{\psi}
dx \ dy \ dt)| \le  C(T)s \lambda  \int_{-T}^T  \int_{\Omega} c  \varphi | \nabla \beta \cdot \nabla \psi|^2 dx \ dy \ dt}}$\\
 $${\displaystyle{ + C(T)s  \lambda  \int_{-T}^T  \int_{\Omega} c  \varphi^3 | \psi|^2 dx \ dy
 \ dt, }}$$
\end{itemize}
where $C(\Omega, \Gamma, T, R_2)$ is a positive constant depending upon
$\Omega, \Gamma, T, R_2$ and $C(T)$ is a positive constant depending upon $T.$
Therefore we obtain the following estimation for $X$:
$$ |X| \le C(\Omega, \Gamma, T, R_2)\left[ (s\lambda^4+s^3\lambda^3)
 \int_{-T}^T  \int_{\Omega}  \varphi^3 |\psi|^2 dx \ dy \ dt + 
s \lambda \int_{-T}^T  \int_{\Omega}  \varphi | \nabla \beta\cdot \nabla \psi|^2 dx \ dy \ dt \right].
$$
The two terms of the previous estimate of $|X|$ are neglectable with respect to
$$ s^3 \lambda^4 \int_{-T}^T  \int_{\Omega}
c^2 \varphi^3 |\nabla \beta|^4 |\psi|^2 dx \ dy \ dt \  
\mbox{ or }\
 s \lambda^2  \int_{-T}^T  \int_{\Omega}
\varphi c^2 | \nabla \beta \cdot \nabla \psi|^2 dx  \  dy \ dt,$$
for $s$ and $\lambda$ sufficiently large. Using Assumption \ref{funct-beta}, we have
$$ 4 s \lambda  \Re \left(\int_{-T}^T  \int_{ \Omega} 
\varphi D^2 \beta(\nabla \psi, \nabla \overline{\psi}) dx \ dy  \ dt\right) 
-2s \lambda \int_{-T}^T  \int_{\Omega}
\varphi\ c \ \nabla c \cdot \nabla \beta |\nabla \psi|^2 dx  \  dy \ dt $$
$$+4 s \lambda^2  \int_{-T}^T  \int_{\Omega}
\varphi c^2 | \nabla \beta \cdot \nabla \psi|^2 dx  \  dy \ dt
\ge C_{pc} s \lambda  \int_{-T}^T  \int_{\Omega}
\varphi  |\nabla \psi|^2 dx  \  dy \ dt,$$
so (\ref{ine}) becomes
$$\|M_1\psi\|^2_{L^2(\widetilde{Q})}  +  \|M_2\psi\|^2_{L^2(\widetilde{Q})} 
 +  4 s^3 \lambda^4 \int_{-T}^T  \int_{\Omega}
c^2 \varphi^3 |\nabla \beta|^4 |\psi|^2 dx \ dy \ dt$$ 
$$+ s \lambda  \int_{-T}^T  \int_{\Omega}
\varphi  |\nabla \psi|^2 dx  \  dy \ dt
\le   2s \lambda \int_{-T}^T  \int_{\partial \Omega} c^2 \varphi \partial_{\nu} \beta
|\partial_{\nu} \psi |^2 \ d\sigma \ dt 
+  \|M\psi\|^2_{L^2(\widetilde{Q})}.$$
Recall that $\partial_{\nu} \beta \leq 0 $ on $ \Gamma^-$, $c(x,y) \in \Lambda(R_1) \cap \Lambda(R_2)$, 
$|\nabla \beta| \in \Lambda(R_1)$ and $\psi=e^{-s \eta} q$, 
then the proof is complete.

\end{proof}
%
\section{Inverse Problem}
%
\setcounter{equation}{0}
In this section, we establish a stability
inequality and deduce a uniqueness result for the
coefficient $c$. 
The Carleman estimate (\ref{Carl}) proved in section $2$ will be the key
ingredient in the proof of such a  stability estimate.\\
Let $q$ be solution of
\begin{equation} \label{syst-q}
\left \{ \begin{array}{lll}
  i \partial_t q +\nabla \cdot(c \nabla q)=0 & \mbox{in} & \Omega \times (0,T),\\
   q(x,y,t)=b(x,y,t)  & \mbox{on} & \partial \Omega \times (0,T),\\
   q(x,y,0)=q_0(x,y) & \mbox{in} & \Omega,
\end{array}\right.
\end{equation}
and $\widetilde{q}$ be solution of
\begin{equation} \label{syst-qtild}
\left \{ \begin{array}{lll}
  i \partial_t \widetilde{q} +\nabla \cdot(\widetilde{c} \nabla \widetilde{q})  =0 & \mbox{in} & \Omega \times (0,T),\\
   \widetilde{q}(x,y,t)=b(x,y,t)  & \mbox{on} & \partial \Omega \times (0,T),\\
   \widetilde{q}(x,y,0)=q_0(x,y) & \mbox{in} & \Omega,
\end{array}\right.
\end{equation}
where $c$ and $\widetilde{c}$ both satisfy Assumption \ref{c}.
If we set $u=q- \widetilde{q}$, $v=\partial_t u$ and $\gamma=\widetilde{c}-c$, then u and v satisfy respectively
\begin{equation} \label{syst-u}
\left \{ \begin{array}{lll}
  i \partial_t u +\nabla \cdot(c \nabla u)= \nabla \cdot (\gamma \nabla \widetilde{q}) 
& \mbox{in} & \Omega \times (0,T),\\
   u(x,y,t)=0  & \mbox{on} & \partial \Omega \times (0,T),\\
   u(x,y,0)=0 & \mbox{in} & \Omega,
\end{array}\right.
\end{equation} 
\begin{equation} \label{syst-v}
\left \{ \begin{array}{lll}
   i \partial_t v +\nabla \cdot(c \nabla v)= \nabla \cdot (\gamma \nabla \partial_t \widetilde{q})=f
 & \mbox{in} & \Omega \times (0,T),\\
   v(x,y,t)=0  & \mbox{on} & \partial \Omega \times (0,T),\\
   v(x,y,0)=\frac{1}{i}\nabla \cdot (\gamma \nabla q_0) & \mbox{in} & \Omega.
\end{array}\right.
\end{equation}
\begin{assumption} \label{q0-bis}
$q_0$ is a real valued function in $\mathcal{C}^3(\Omega)$
\end{assumption}
We extend the function $v$ on $\Omega \times (-T,T)$ by the formula
$v(x,y,t)=-\overline{v}(x,y,-t)$ for every $(x,y,t) \in \Omega \times (-T,0)$.
Note that this extension is available if the initial data is a real valued function. For
a pure imaginary initial data, the right extension is $v(x,y,t)=\overline{v}(x,y,-t)$.
Note that these extensions satisfy the previous Carleman estimate.
%
\subsection{Energy Estimate}
%
We assume throughout this section that $\gamma \in H^1_0(\Omega)$.
We introduce
\begin{equation}\label{E+I}
\mathbb{E}(t)=\int_{\Omega}   e^{-2s \eta(x,y,t)} |\partial_t u(x,y,t)|^2 dx \ dy
+\int_{\Omega} \varphi^{-1}(x,y,t) \  e^{-2s \eta(x,y,t)} |\partial_t \nabla u(x,y,t)|^2 dx \ dy.
\end{equation}
In this section, we will give an estimation of $\mathbb{E}(0)$.\\ \\ \noindent
%
%
{\bf First Step}: We first give an estimation of $\int_{\Omega}   e^{-2s \eta(x,y,0)} |\partial_t u(x,y,0)|^2 dx \ dy$.\\ \noindent
We set $\psi = e^{-s \eta} v$. 
With the operator
\begin{equation}
\label{eq: M1} 
 M_1\psi  = i\partial_t \psi+\nabla \cdot (c \nabla\psi)
+s^2 |\nabla \eta |^2 \psi,
\end{equation}
we introduce, following \cite{BP}, 
\begin{eqnarray*}
  \mathcal{I} = 2\Im \left(\int_{-T}^{0} \int_{\Omega} M_1 \psi\;\overline{\psi}\;dx \ dy \ dt\right).
\end{eqnarray*}
\begin{assumption} \label{regu} 
$\partial_t \widetilde{q}, \ \nabla(\partial_t \widetilde{q}), \
\Delta(\partial_t \widetilde{q}) $ are in $\Lambda(R_2)$.
\end{assumption}

We have the following estimate
\begin{lemma}
  \label{lemma1}
 We assume that  Assumption \ref{regu} is satisfied. 
 Then there exists a positive constant $C=C(\Omega, \Gamma, T, R_1, R_2)$ such that
 for any $\lambda \geq \lambda_0$ and $s \geq s_0$, we have
  \begin{eqnarray*} \label{I2}
  \mathcal{I} & = & \int_{\Omega}  e^{-2s \eta(x,y,0)} |\partial_t u(x,y,0)|^2 \ d x \ dy\\ 
  \mbox{ and} \\ \nonumber
    |\mathcal{I}| & \leq & C s^{-3/2} \lambda^{-2}
   \int_{\Omega}  e^{-2s \eta(x,y,0)}( |\gamma|^2 +|\nabla \gamma|^2)\ dx \ dy \\ 
     & & + C s^{-1/2} \lambda^{-1} \int_{-T}^T \int_{\Gamma^+}  e^{-2s \eta}\ \varphi\ 
\partial_{\nu} \beta \ |\partial_{\nu} v|^2\ d\sigma \ dt.
  \end{eqnarray*}
\end{lemma}
\begin{proof}
In a first step, we calculate $\mathcal{I}$
\begin{eqnarray*}
\mathcal{I} & = & 2\Im \left(\int_{-T}^{0} \int_{\Omega}\left( i \partial_t \psi\ \overline{\psi}
+ \nabla \cdot(c\nabla \psi)\ \overline{\psi}+ s^2|\nabla \eta|^2 \psi\ \overline{\psi}\right)
\ d x \ dy \ d t \right)\\
& = & 2\Re \left(\int_{-T}^{0} \int_{\Omega}  \partial_t\ \psi \overline{\psi} \ d x \ dy \ d t\right)\\
& = &  \int_{-T}^{0} \int_{\Omega}  \partial_t |\psi|^2 \ d x \ dy \ d t\\
& = &  \int_{\Omega}|\psi(x,y,0)|^2 \ d x \ dy \\
& = &  \int_{\Omega}  e^{-2s \eta(x,y,0)} |v(x,y,0)|^2 \ d x \ dy.
\end{eqnarray*}
So, we have
$$\mathcal{I} =  \int_{\Omega}  e^{-2s \eta(x,y,0)} |\partial_t u(x,y,0)|^2 \ d x \ dy.$$
In a second step, we estimate $\mathcal{I}$. Using Young inequality we
can write
\begin{eqnarray*}
|\mathcal{I}| & \leq &  2\left( \int_{-T}^{T} \int_{\Omega} 
|M_1 \psi|^2 \ d x \ dy \ d t\right)^{\frac{1}{2}}
\left( \int_{-T}^{T} \int_{\Omega} 
e^{-2s \eta}|v|^2 \ d x \ dy \ d t\right)^{\frac{1}{2}}\\
& \leq & C(T) s^{-\frac{3}{2}} \lambda^{-2} \left(
\|M_1 \psi \|^2_{L^2(\widetilde{Q})} +
s^3 \lambda^4  \int_{-T}^{T} \int_{\Omega} 
e^{-2s \eta} \varphi^3 |v|^2 \ d x \ dy \ d t\right)
\end{eqnarray*}
with $C(T)$ a positive constant which depends on $T.$
Then with the Carleman estimate (\ref{Carl}) proved in section $2$ we have
\begin{eqnarray*}
|\mathcal{I}| & \leq & C  s^{-\frac{3}{2}} \lambda^{-2} \left(
s \lambda  \int_{-T}^{T} \int_{\Gamma^+} 
e^{-2s \eta} \varphi \partial_{\nu} \beta |\partial_n v|^2 \ d \sigma \ d t
+\int_{-T}^{T} \int_{\Omega} e^{-2s \eta}
|\nabla \cdot (\gamma \nabla \partial_t \widetilde{q}) \ d x \ dy \ d t \right),
\end{eqnarray*}
where $C=C(\Omega, \Gamma^+, T, R_1, R_2)$ is a positive constant.
Using Assumption \ref{regu}, since 
$$e^{{-2s \eta}(x,y,t)} \le e^{{-2s \eta}(x,y,0)} \;\;\mbox{ for all} \; t \in (-T,T),$$
we obtain for $s$ and $\lambda$ sufficiently large the estimate
\begin{eqnarray} \label{est1-I}
|\mathcal{I}| \leq C\ s^{-3/2} \lambda^{-2}
\int_{\Omega} e^{{-2s \eta}(x,y,0)} 
(| \gamma |^2+|\nabla \gamma|^2) \ d x \ dy \\ \nonumber
+C\ s^{-1/2} \lambda^{-1} \int_{-T}^T \int_{\Gamma^+}  e^{-2s \eta} \varphi\ \partial_{\nu} \beta \ |\partial_{\nu} v |^2\ d\sigma \ dt.
\end{eqnarray}
where $C=C(\Omega, \Gamma^+, T, R_1, R_2)$ is a positive constant.
\end{proof}
%
%

 \noindent
{\bf Second Step}:  We then give an estimate of $\int_{\Omega} \varphi^{-1}(x,y,0) \   e^{-2s \eta(x,y,0)} |\partial_t \nabla u(x,y,0)|^2 dx \ dy$.\\ \noindent
We denote
\begin{equation} \label{en}
E(t):=\int_{\Omega}  c\  \varphi^{-1}(x,y,t) \ e^{-2s \eta(x,y,t)} |\nabla v(x,y,t)|^2 dx \ dy,
\end{equation}
where $\varphi^{-1}=\frac{1}{\varphi}$.
We give an estimate for $E(0)$ in Theorem \ref{th-EE}.
In a first step we prove the following lemma :

\begin{lemma} \label{lem4}
Let $v$ be solution of (\ref{syst-v}) in the following class
$$v \in C([0,T],H^1(\Omega)),\ \partial_{\nu} v \in L^2(0,T,L^2(\Gamma)).$$
Then the following identity holds true
$$\ E(\tau)-E(\kappa)=-2 \Re\left( \int_{\kappa}^{\tau} \int_{\Omega} 
e^{-2s \eta}f\ \varphi^{-1} \partial_t \overline{v} \ dx\ dy\ dt\right)$$
$$+\Re\left(\int_{\kappa}^{\tau} \int_{\Omega}c\
  e^{-2s\eta} (-4s\lambda +2\lambda \varphi^{-1})
\partial_t \overline{v}\ \nabla \beta \cdot \nabla v \ dx\ dy\ dt \right)$$
$$-2s\int_{\kappa}^{\tau} \int_{\Omega}c\ e^{-2s\eta}\varphi^{-1} 
\partial_t \eta | \nabla v|^2 \ dx\ dy\ dt  +\int_{\kappa}^{\tau}
\int_{\Omega} ce^{-2s\eta} \partial_t (\varphi^{-1}) |\nabla v|^2,$$
for $f \in H^1_0(\Omega)$.
\end{lemma}
\begin{proof}
Since $v$ is solution of $(3.4)$ note that $\partial_t
\overline{v}=i\overline{f} -i\nabla \cdot (c\nabla \overline{v}).$
Therefore, we obtain the two following equalities. 
\begin{equation} \label{ii}
\ \Re\int_{\kappa}^{\tau} \int_{\Omega} e^{-2s \eta} f\ 
\varphi^{-1} 
\partial_t\overline{v} \ dx\ dy\ dt=
\Re\left(-i\int_{\kappa}^{\tau} \int_{\Omega} e^{-2s \eta}f\ \varphi^{-1} 
\nabla \cdot (c\ \nabla \overline{v}) \ dx\ dy\ dt\right),
\end{equation}
\begin{eqnarray} \label{iii}
\Re\left(\int_{\kappa}^{\tau} \int_{\Omega} c\ e^{-2s \eta}
\partial_t \overline{v}\ (-4s\lambda +2\lambda \varphi^{-1}) 
\nabla \beta \cdot 
\nabla v \ dx\ dy\ dt\right)\\ \nonumber
= \Re\left(i \int_{\kappa}^{\tau} \int_{\Omega} c\ e^{-2s \eta}\overline{f}\
  (-4s\lambda +2\lambda \varphi^{-1}) \nabla \beta \cdot \nabla v\ dx\ dy\ dt\right)\\ \nonumber
-\Re\left(i \int_{\kappa}^{\tau} \int_{\Omega}c\ e^{-2s
    \eta} (-4s\lambda +2\lambda \varphi^{-1}) 
\nabla\cdot(c \nabla \overline{v}) \nabla \beta \cdot \nabla v \ dx\ dy\ dt\right).
\end{eqnarray}
 We multiply the first equation of (\ref{syst-v}) by $e^{-2s \eta}\
 \varphi^{-1} \partial_t \overline{v}$ and
we integrate on $(\kappa,\tau)\times \Omega$, where $[\kappa,\tau] \subset [0,T]$. So, if we consider the real 
part of the obtained equality, we have
$$0=\Re\left( i\int_{\kappa}^{\tau} \int_{\Omega}e^{-2s \eta}\ \varphi^{-1}
|\partial_t v|^2\ dx\ dy\ dt\right)
=\Re\left(\int_{\kappa}^{\tau} \int_{\Omega}e^{-2s \eta}\ f\ \varphi^{-1} \partial_t
  \overline{v}\ dx\ dy\ dt\right)$$
$$-\Re\left(\int_{\kappa}^{\tau} \int_{\Omega}e^{-2s \eta}\ \varphi^{-1} 
\nabla \cdot(c\nabla v)\partial_t \overline{v}\ dx\ dy\ dt\right).$$
Then by integration by parts, we obtain
$$0 =  \Re \left(\int_{\kappa}^{\tau} \int_{\Omega}e^{-2s \eta}\ f\
  \varphi^{-1} \partial_t \overline{v}\ dx\ dy\ dt\right)
+\Re\left( \int_{\kappa}^{\tau} \int_{\Omega} \ c \nabla v \cdot \nabla 
(e^{-2s \eta}\ \varphi^{-1} \partial_t \overline{v})\ dx\ dy\ dt\right)$$
$$0 = \Re \left( \int_{\kappa}^{\tau} \int_{\Omega}e^{-2s \eta}\ f\
  \varphi^{-1} \partial_t \overline{v}\ dx\ dy\ dt\right)
-2s \Re \left( \int_{\kappa}^{\tau} \int_{\Omega} \ c\ \varphi^{-1} e^{-2s \eta}\
\partial_t \overline{v} \nabla v \cdot \nabla \eta \ dx\ dy\ dt \right) $$
$$-\lambda \Re \left(\int_{\kappa}^{\tau} \int_{\Omega}c\ e^{-2s \eta}\ 
\varphi^{-1} \partial_t \overline{v} \nabla v \cdot \nabla \beta 
\ dx\ dy\ dt\right)
+\Re \left(\int_{\kappa}^{\tau} \int_{\Omega}c\ e^{-2s \eta}\ 
\varphi^{-1} \nabla v \cdot \nabla (\partial_t \overline{v})
\ dx\ dy\ dt\right).$$
Note that
$$E(\tau)-E(\kappa)=\int_{\Omega} \int_{\kappa}^{\tau} \ c\ \partial_t (e^{-2s
  \eta}\  \varphi^{-1} |\nabla v|^2)\ dx\ dy\ dt.$$
Therefore we have
$$0 = \Re \left( \int_{\kappa}^{\tau} \int_{\Omega}e^{-2s \eta}\ f\
  \varphi^{-1} \partial_t \overline{v}\ dx\ dy\ dt\right)
-2s \Re \left( \int_{\kappa}^{\tau} \int_{\Omega} \ c\ \varphi^{-1} e^{-2s \eta}\
\partial_t \overline{v} \nabla v \cdot \nabla \eta \ dx\ dy\ dt \right) $$
$$-\lambda \Re \left(\int_{\kappa}^{\tau} \int_{\Omega}c\ e^{-2s \eta}\ 
\varphi^{-1} \partial_t \overline{v} \nabla v \cdot \nabla \beta 
\ dx\ dy\ dt\right)
+\frac{1}{2} E(\tau) -\frac{1}{2} E(\kappa)$$
$$+s \int_{\kappa}^{\tau} \int_{\Omega}c\ e^{-2s \eta}\ 
\varphi^{-1} \partial_t \eta |\nabla v|^2 \ dx\ dy\ dt
-\frac{1}{2} \int_{\kappa}^{\tau} c e^{-2s\eta} \partial_t (\varphi^{-1})
|\nabla v|^2 \ dx \ dy \ dt,$$
and the proof of Lemma \ref{lem4} is complete.
\end{proof}
%
%
\begin{theorem}
\label{th-EE}
Let $v$ be solution of (\ref{syst-v}) in the following class
$$v \in C([0,T],H^1(\Omega)),\ \partial_{\nu} v \in L^2(0,T,L^2(\Gamma)).$$
We assume that Assumptions \ref{c} and \ref{funct-beta} are checked.  
Then there exists a positive constant \\$C=C(\Omega, \Gamma,T, R_1, R_2) >0$ such that
\begin{equation} \label{EE}
E(0) \le C \left[ s^2 \lambda^2 
\int_{-T}^T \int_{\Gamma^+}e^{-2s \eta} \varphi\ \partial_{\nu}\beta  
\ |\partial_{\nu} v|^2\ d\sigma \ dt+s\lambda \int \int_Q e^{-2s\eta} |f|^2 
\right] 
\end{equation}
for $s$ and $\lambda$ sufficiently large.
\end{theorem}
\begin{proof}
We apply Lemma \ref{lem4} with $\kappa=0$ and $\tau=T$. Since $E(T)=0$, we obtain
\begin{eqnarray*}
E(0)=2 \Re\left( \int_{0}^{T} \int_{\Omega} 
e^{-2s \eta}f\ \varphi^{-1} \partial_t \overline{v} \ dx\ dy\ dt\right)
-\Re\left(\int_{0}^{T} \int_{\Omega}c\
  e^{-2s\eta} (-4s\lambda +2\lambda \varphi^{-1})
\partial_t \overline{v}\ \nabla \beta \cdot \nabla v \ dx\ dy\ dt \right)\\
+2s\int_{0}^{T} \int_{\Omega}c\ e^{-2s\eta}\varphi^{-1} 
\partial_t \eta | \nabla v|^2 \ dx\ dy\ dt  
-\int_{0}^{T}
\int_{\Omega} ce^{-2s\eta} \partial_t (\varphi^{-1}) |\nabla v|^2.
\end{eqnarray*}
We give now estimates of the four integrals in the previous equality.\\ \\ \noindent
{\bf First integral}:  $B:=2 \Re\left( \int_{0}^{T} \int_{\Omega} 
e^{-2s\eta} \varphi^{-1} f\ \partial_t \overline{v} \ dx\ dy\ dt \right)$.\\ \noindent
Using (\ref{ii}), we have :
\begin{equation} \label{est-f}
B=2 \Re\left( \int_{0}^{T} \int_{\Omega} e^{-2s\eta} \varphi^{-1} 
f\ \partial_t \overline{v} \ dx\ dy\ dt \right)=
\Re\left( \int_{0}^{T} \int_{\Omega} e^{-2s\eta} \varphi^{-1} 
f\ \partial_t \overline{v} \ dx\ dy\ dt \right)
\end{equation}
$$+\Re\left(-i\int_{0}^{T} \int_{\Omega} e^{-2s\eta} \varphi^{-1} 
f\ \nabla \cdot (c\ \nabla \overline{v}) \ dx\ dy\ dt\right).$$
Recall that if we set $\psi=e^{-s \eta}v$, then $M \psi = e^{-s \eta} H(e^{s \eta} \psi)=M_1 \psi + M_2 \psi$ for $s>0$
with
$$
  M_1\psi  : = i\partial_t \psi+\nabla \cdot (c \nabla\psi)
+s^2 c |\nabla \eta |^2 \psi,
$$
$$
  M_2\psi  : = is \partial_t\eta \psi+2 c s \nabla \eta  \cdot \nabla \psi
+s \nabla \cdot (c \nabla \eta) \psi.
$$
So (\ref{est-f}) becomes
$$
B=2 \Re\left( \int_{0}^{T} \int_{\Omega} e^{-2s\eta} \varphi^{-1} 
f\ \partial_t \overline{v} \ dx\ dy\ dt \right)=
 \Re\left( \int_{0}^{T} \int_{\Omega} e^{-2s\eta} \varphi^{-1} 
f\ e^{s \eta}(s \partial_t \eta \overline{\psi}
+\partial_t \overline{\psi}) \ dx\ dy\ dt \right)$$
$$+\Re\left( -i \int_{0}^{T} \int_{\Omega} e^{-2s\eta} \varphi^{-1} 
f\ e^{s \eta}(s\nabla \cdot(c\nabla \eta)\overline{\psi}
+c\ s^2|\nabla \eta|^2 \overline{\psi} +2c\ s\ \nabla \eta \cdot \nabla \overline{\psi}
+\nabla \cdot(c\nabla \overline{\psi}))\ dx\ dy\ dt \right)$$
$$=\Re\left( -i \int_{0}^{T} \int_{\Omega} e^{-2s\eta} \varphi^{-1} 
f\ e^{s \eta}(i \partial_t \overline{\psi}+
\nabla \cdot ( c\nabla \overline{\psi} ))\ dx\ dy\ dt \right)$$
$$+\Re\left( \int_{0}^{T} \int_{\Omega} e^{-2s\eta} \varphi^{-1} f\ e^{s \eta}
(s \partial_t \eta -i\ s\nabla \cdot(c\nabla \eta)
-i\ c\ s^2|\nabla \eta|^2 ) \overline{\psi}\ dx\ dy\ dt \right)$$
$$+\Re\left( -i \int_{0}^{T} \int_{\Omega} e^{-2s\eta} \varphi^{-1} f\ e^{s \eta}
2c\ s\ \nabla \eta \cdot \nabla \overline{\psi}\ dx\ dy\ dt \right)
$$
Note that
$$i\partial_t \overline{\psi}+\nabla \cdot (c \nabla\overline{\psi})
= M_1 \overline{\psi}  -s^2 c |\nabla \eta |^2 \overline{\psi}.$$
Then, we obtain
$$B
=\Re\left( -i \int_{0}^{T} \int_{\Omega} e^{-2s\eta} \varphi^{-1} f\ e^{s \eta}
(M_1 \overline{\psi}  -s^2 c |\nabla \eta |^2 \overline{\psi})\ dx\ dy\ dt \right)$$
$$+\Re\left( \int_{0}^{T} \int_{\Omega} e^{-2s\eta} \varphi^{-1} f\ e^{s \eta}
(s \partial_t \eta -i\ s\nabla \cdot(c\nabla \eta)
-i\ c\ s^2|\nabla \eta|^2 ) \overline{\psi}\ dx\ dy\ dt \right)$$
$$+\Re\left( -i \int_{0}^{T} \int_{\Omega} e^{-2s\eta} \varphi^{-1} f\ e^{s \eta}
2c\ s\ \nabla \eta \cdot \nabla \overline{\psi}\ dx\ dy\ dt \right).$$
If we come back to the function $v$, the previous equality becomes :
$$B
=\Re\left( -i \int_{0}^{T} \int_{\Omega} (f
e^{-s \eta} \varphi^{-1} M_1 (e^{-s \eta}\overline{v})  
-s^2 c |\nabla \eta |^2 e^{-2s\eta} \varphi^{-1} f\ \overline{v})\ dx\ dy\ dt \right)$$
$$+\Re\left( \int_{0}^{T} \int_{\Omega} e^{-2s\eta} \varphi^{-1} f\ 
(s \partial_t \eta -i\ s\nabla \cdot(c\nabla \eta)
+i\ c\ s^2|\nabla \eta|^2 ) \overline{v}\ dx\ dy\ dt \right)$$
$$+\Re\left( -i \int_{0}^{T} \int_{\Omega} e^{-2s\eta} \varphi^{-1} f\ 
2c\ s\ \nabla \eta \cdot \nabla \overline{v}\ dx\ dy\ dt \right).$$
Then there exists a positive constant $C=C(\Omega, \Gamma, T, R_1, R_2)$ such that:
\begin{equation}
|\int_{0}^{T} \int_{\Omega} e^{-2s\eta} \varphi^{-1} 
f\ \partial_t \overline{v} \ dx\ dy\ dt|
\le C \left[ s\lambda \int \int_{Q} e^{-2s\eta} 
|f|^2\ dx\ dy\ dt
+ \|M_1(e^{-s \eta}\overline{v})\|^2_{L^2(Q)}\right. 
\end{equation}
$$ \left.+s^3 \lambda^4 \int \int_{Q}e^{-2 s \eta}\ \varphi^2 |v|^2 \ dx\ dy\ dt 
+s\lambda \int \int_{Q}  e^{-2s\eta} 
|\nabla v|^2 \ dx\ dy\ dt \right]. $$
\\ \noindent
{\bf Second integral}:
$D:=\Re\left(\int_{0}^{T} \int_{\Omega}c\
  (-4s\lambda +2\lambda \varphi^{-1}) \ e^{-2s\eta} (1+\varphi^{-1})
\partial_t \overline{v}\ \nabla \beta \cdot \nabla v \ dx\ dy\ dt
\right).$ \\ \noindent
We denote by $\rho:=-4s\lambda+ 2\lambda \varphi^{-1}.$ Using (\ref{iii}), we
have 
$$2D=\Re\left(\int_{0}^{T} \int_{\Omega}c\
  e^{-2s\eta} \rho 
\partial_t \overline{v}\ \nabla \beta \cdot \nabla v \ dx\ dy\ dt
\right)
+ 
\Re\left(i \int_{0}^{T} \int_{\Omega} c\ e^{-2s \eta}\overline{f}\
  \rho \nabla \beta \cdot \nabla v\ dx\ dy\ dt\right)$$
$$- \Re\left(i \int_{0}^{T} \int_{\Omega}c\ e^{-2s
    \eta}\rho
\nabla\cdot(c \nabla \overline{v}) \nabla \beta \cdot \nabla v \ dx\ dy\ dt\right).$$
If we introduce $\psi=e^{-s\eta}v,$ we get: 
$$2D=  \Re \left(
  \int_{0}^{T} \int_{\Omega} c e^{-2s\eta} \rho \nabla v \cdot
  \nabla \beta (s\partial_t \eta e^{s\eta} \overline{\psi} +e^{s\eta}
  \partial_t \overline{\psi}) \ dx\ dy \ dt \right)
+ \Re \left( i\int_{0}^{T} \int_{\Omega}
  ce^{-2s\eta} \nabla v \cdot \nabla \beta \ \overline{f} \ dx \ dy \ dt
\right) $$
$$-\Re \left( i\int_{0}^{T} \int_{\Omega} 
ce^{-2s\eta} \rho \nabla v \cdot \nabla \beta [s^2c \overline{\psi} |\nabla
\eta|^2 e^{s\eta} +2sce^{s\eta} \nabla \overline{\psi} \cdot \nabla \eta +s
e^{s\eta} \overline{\psi} \nabla \cdot (c\nabla \eta) +e^{s\eta} \nabla
\cdot (c\nabla \overline{\psi})]\right).$$
Therefore
$$2D= \Re \left( -i\int_{0}^{T} \int_{\Omega}
  ce^{-s\eta} \rho \nabla v\cdot \nabla \beta [ M_1 \overline{\psi} +M_2
  \overline{\psi}] \ dx \ dy \ dt \right)+ \Re \left( i\int_{0}^{T} \int_{\Omega}
  ce^{-2s\eta} \rho \nabla v\cdot \nabla \beta \  \overline{f} \ dx \ dy \
  dt \right).$$
Thus there exists a positive constant $C=C(\Omega, \Gamma,T, R_1, R_2)$ such that 
\begin{equation}
|2D|\leq C \left[ s\lambda \int \int_Q [|M_1(e^{-s\eta} \overline{v})|^2 +
|M_2(e^{-s\eta} \overline{v})|^2] \ dx \ dy \ dt \right.
\end{equation}
$$ \left. +s\lambda \int \int_Q e^{-2s\eta}  |\nabla v|^2   
+s\lambda \int \int_Q
e^{-2s\eta} |f|^2 \ dx \ dy \ dt \right].$$
\\ \noindent
{\bf Two last integrals}: \\ \noindent
There exists a positive constant $C=C(\Omega, \Gamma,T, R_2)$ such that
\begin{equation}
|\int_{0}^{T} \int_{\Omega} c e^{-2s\eta} \partial_t
(\varphi^{-1}) |\nabla v|^2  \ dx \ dy \ dt | \leq C\int \int_Q e^{-2s\eta}
|\nabla v|^2 \ dx \ dy \ dt, 
\end{equation}
and  
\begin{equation}
|s\int_{0}^{T} \int_{\Omega} c\partial_t \eta e^{-2s\eta}
\varphi^{-1} |\nabla v|^2 \ dx \ dy \ dt | \leq C\ s \int \int_Q
\varphi e^{-2s\eta} |\nabla v|^2 \ dx \ dy \ dt.
\end{equation}
Using now the Carleman estimate of Theorem \ref{th-Carl} and Lemma \ref{lem4}, 
from (3.9)-(3.12), we deduce the existence of 
a positive constant
$C=C(\Omega, \Gamma,T, R_1, R_2)$ such that:
$$E(0) \le    
 C \left[s^2 \lambda^2 \int_{-T}^{T} \int_{\Gamma^+} e^{-2s \eta} 
\varphi\ \partial_{\nu}\beta\ 
|\partial_{\nu} v|^2\ d\sigma\ dt
+s\lambda \int \int_Q e^{-2s\eta}  | f|^2 \ dx \ dy \ dt \right],$$
and the proof is complete.
\end{proof}
%
\subsection{Stability Estimate}
%
Now following an idea developed in \cite{IIY} for Lam\'e system in bounded domains,
we give an underestimate for $\mathbb{E}(0)$.
We adapt the proof of lemma 3.2 of \cite{IIY} to an unbounded domain.
\begin{assumption} \label{q0}
\begin{itemize}
\item
$q_0$ and all its derivatives up to order three are in $\Lambda(R_2)$
\item
$|\nabla \beta \cdot \nabla q_0|\in \Lambda(R_1)$
\end{itemize}
\end{assumption}
\vskip 0.3cm
%
%
\begin{lemma}\label{P0g}
We consider the first order partial differential operator
$$(P_0 g)(x,y)= \partial_x q_0(x,y) \partial_x g(x,y) 
+ \partial_y(x,y) \partial_y g(x,y), P_0g:=\nabla q_0 \cdot \nabla g$$
where $q_0$ satisfies Assumptions \ref{q0-bis}, \ref{q0}.
Then there exist positive constants $\lambda_1>0$, $s_1>0$ and $C=C(\Omega, \Gamma,T, R_1, R_2)$
such that for all $\lambda \geq \lambda_1$ and $s \geq s_1$
$$s^2  \lambda^2\int_{\Omega}  \varphi_0 \ e^{-2s \eta_0} |g|^2\ dx \ dy
\leq C  \int_{\Omega} \varphi^{-1}_0 \ e^{-2s \eta_0} \ |P_0 g|^2\ dx \ dy$$
with $\eta_0(x,y):=\eta(x,y,0)$, $\varphi_0 (x,y):=\varphi(x,y,0)$ and 
for $g \in H^1_0(\Omega)$.
\end{lemma}

\begin{proof}
Let $g \in H^1_0(\Omega)$. We denote by $w=e^{-s\eta_0}g$ with
$\eta_0:=\eta(x,y,0)$ and 
$Q_0w= e^{-s \eta_0}P_0(e^{s \eta_0}w),$ so we get $Q_0w=sw P_0 \eta_0 +P_0w$.
Therefore, we have:
$$\int_{\Omega}  \varphi_0^{-1} \ Q_0 w\  \overline{Q_0 w}\ dx \ dy =
s^2 \int_{\Omega} \varphi_0^{-1} \  |w|^2|P_0 \eta_0|^2 \ dx \ dy
+\int_{\Omega} \varphi_0^{-1} \ |P_0 w|^2 \ dx \ dy$$
$$+2s \Re \left(\int_{\Omega} \ \varphi_0^{-1} \ 
w\ P_0 \eta_0\ \overline{P_0 w}\ dx \ dy \right)$$
$$=s^2  \lambda^2 \int_{\Omega} \varphi_0 \ |w|^2 (\nabla q_0 \cdot \nabla 
\beta)^2 \ dx \ dy
+\int_{\Omega}  \varphi_0^{-1} \ |P_0 w|^2 \ dx \ dy
-s  \lambda \int_{\Omega} \nabla q_0 \cdot \nabla \beta  \ 
\nabla q_0 \cdot \nabla (|w|^2) \ dx \ dy .$$
So, integrating by parts, we obtain 
$$\int_{\Omega}  \varphi_0^{-1} \ Q_0 w\  \overline{Q_0 w}\ dx \ dy =
s^2  \lambda^2 \int_{\Omega}  \varphi_0 \ |w|^2 (\nabla q_0 \cdot \nabla 
\beta)^2 \ dx \ dy
+\int_{\Omega} \varphi_0^{-1} \ \ |P_0 w|^2 \ dx \ dy$$
$$+s \lambda \int_{\Omega} \nabla \cdot (\nabla q_0 \cdot \nabla \beta 
\nabla q_0 )|w|^2 \ dx \ dy .$$
Thus
$$\int_{\Omega} \varphi_0^{-1}  \ e^{-2s \eta_0} |P_0 g|^2 \ dx \ dy
\geq s^2 \lambda^2\int_{\Omega}  \varphi_0 \ 
|\nabla \beta \cdot \nabla q_0|^2e^{-2s \eta_0} |g|^2\ dx \ dy
+s  \lambda \int_{\Omega} \nabla \cdot(P_0 \beta\ \nabla q_0) e^{-2s\eta_0}|g|^2\ dx \ dy.$$
Using Assumptions \ref{q0-bis} and \ref{q0}, we can conclude for $s$ and $\lambda$
 sufficiently large.
\end{proof}
Then, we deduce the following result.
\begin{lemma}\label{Carl3}
Let $u$ be solution of (\ref{syst-v}). We assume that Assumptions 
\ref{funct-beta}, \ref{q0-bis} and \ref{q0}
are satisfied. Then there exists a positive constant
$C=C(\Omega, \Gamma,T, R_1, R_2)$ such that for $s$ and $\lambda$ sufficiently large, the two following estimates 
hold true
\begin{equation} \label{est1-gam}
s^2 \lambda^2 \int_{\Omega} \varphi_0 \ e^{-2s \eta_0} |\gamma|^2\ dx \ dy
\leq C  \int_{\Omega} \varphi_0^{-1}  \ e^{-2s \eta_0} 
|\partial_t u(x,y,0)|^2\ dx \ dy,
\end{equation}
\begin{equation} \label{est2-gam}
s^2  \lambda^2\int_{\Omega} e^{-2s \eta_0}|\nabla \gamma|^2\ dx \ dy
\leq C  \int_{\Omega} \varphi_0^{-1} \ e^{-2s \eta_0} \left(|\nabla(\partial_t u(x,y,0))|^2+|\gamma|^2\right)\ dx \ dy,
\end{equation}
for $\gamma \in H^2_0(\Omega)$.
\end{lemma}
\begin{proof}
We apply Lemma \ref{P0g} to the first order partial differential equations satisfied by
\begin{itemize}
\item
$\gamma$ 
given by the initial condition in (\ref{syst-v})
$$
P_0 \gamma:=\partial_x q_0 \partial_x \gamma +\partial_y q_0 \partial_y \gamma
=i \partial_t u(x,y,0)-\gamma \Delta q_0,
$$
\item
$\partial_x \gamma$
given by the $x$-derivative of the initial condition in (\ref{syst-v})
\begin{eqnarray*} \label{ci_x}
P_0 \partial_x \gamma & := & \partial_x q_0 \partial_x(\partial_x \gamma) +\partial_y q_0 \partial_y (\partial_x \gamma)\\ 
& = & i \partial_t (\partial_x u(x,y,0))-\partial_x \gamma (\Delta q_0+\partial_{xx}q_0)
-\partial_y \gamma \partial_{xy} q_0
-\gamma \partial_x (\Delta q_0),
\end{eqnarray*}
\item
$\partial_y \gamma$
given by the $y$-derivative of the initial condition in (\ref{syst-v})
\begin{eqnarray*} \label{ci_y}
P_0 \partial_y \gamma & := & \partial_x q_0 \partial_x(\partial_y \gamma) +\partial_y q_0 \partial_y (\partial_y \gamma)\\
& = & i \partial_t (\partial_y u(x,y,0))-\partial_y \gamma (\Delta q_0+\partial_{yy}q_0)
-\partial_x \gamma \partial_{xy} q_0
-\gamma \partial_y (\Delta q_0).
\end{eqnarray*}
\end{itemize}
Then using Lemma \ref{P0g} and Assumptions \ref{q0-bis}, \ref{q0}, the proof of Lemma \ref{Carl3} is complete.
\end{proof}
\begin{theorem} \label{stab}
 Let $q$ and $\widetilde{q}$ be solutions of (\ref{syst-q}) and (\ref{syst-qtild})
such that $c-\widetilde{c} \in H^2_0(\Omega)$.
We assume that  Assumptions \ref{c}, \ref{funct-beta}, \ref{regu}, \ref{q0-bis} and \ref{q0} are satisfied.
  Then there exists a positive constant
  $C=C(\Omega, \Gamma,T, R_1, R_2)$ 
such that for $s$ and $\lambda$ sufficiently large,
$$\int_{\Omega} \varphi_0 \ e^{-2s\eta_0} (|c-\widetilde{c}|^2 +|\nabla (c-\widetilde{c}) |^2) \
dx \ dy 
\leq C  \int_{-T}^T \int_{\Gamma^+} \varphi \ e^{-2s\eta} \partial_{\nu}
\beta \ |\partial_{\nu} (\partial_t q -\partial_t \widetilde{q})|^2 \ d\sigma \ dt.$$
 \end{theorem}
\begin{proof}
Adding (\ref{est1-gam}) and (\ref{est2-gam}) we obtain 
using the estimate (\ref{est1-I}) for  $|\mathcal{I}|$ and
the energy estimate (\ref{EE}) for $E(0)$
$$
s^2 \lambda^2\int_{\Omega} \varphi_0 \ e^{-2s \eta_0}
\left(|\nabla \gamma|^2+|\gamma|^2\right)\ dx \ dy
  \leq  C  \int_{\Omega} \varphi_0^{-1} \ e^{-2s \eta_0} 
\left(|\nabla(\partial_t u(x,y,0))|^2+|\partial_t u(x,y,0)|^2\right)
\ dx \ dy$$
$$ \leq  C(|\mathcal{I}| +E(0)) $$
$$\leq 
C s^{-3/2} \lambda^{-2}
\int_{\Omega} e^{-2s \eta_0} 
(| \gamma |^2+|\nabla \gamma|^2) \ d x \ dy 
 +  C s^{-1/2} \lambda^{-1} \int_{-T}^T \int_{\Gamma^+}  e^{-2s \eta} \varphi\ 
\partial_{\nu} \beta \ |\partial_{\nu} v |^2\ d\sigma \ dt$$
$$+C  s^2 \lambda^2 
\int_{-T}^T \int_{\Gamma^+}e^{-2s \eta} \varphi\ \partial_{\nu}\beta  
\ |\partial_{\nu} v|^2\ d\sigma \ dt  +  C s\lambda \int \int_Q e^{-2s\eta} |f|^2.$$
So we get
$$
s^2 \lambda^2 \int_{\Omega} \varphi_0 \ e^{-2s \eta_0} 
\left(|\nabla \gamma|^2+|\gamma|^2\right)\ dx \ dy
 \leq  C s^2 \lambda^{2} \int_{-T}^T \int_{\Gamma^+}  e^{-2s \eta}\ \varphi\ 
\partial_{\nu} \beta \ |\partial_{\nu} v|^2\ d\sigma \ dt$$
$$+C s \lambda  \int  \int_{Q}  e^{-2s \eta}\ |\nabla \cdot (\gamma \nabla \partial_t
\widetilde{q} )|^2 \ dx \ dy \ dt$$
$$\leq  C s^2 \lambda^{2} \int_{-T}^T \int_{\Gamma^+}  e^{-2s \eta}\ \varphi\ 
\partial_{\nu} \beta \ |\partial_{\nu} v|^2\ d\sigma \ dt
+C s\lambda \int  \int_{Q}  e^{-2s \eta}\ \left( |\nabla \gamma |^2 + 
| \gamma|^2 \right)
\ dx \ dy \ dt.$$
Then, for $s$ and $\lambda$ sufficiently large, the theorem is proved.
\end{proof}
\begin{remark}
This result is also available for the heat equation in bounded or unbounded domains.\\
Note that all the previous results proved in $\mathbb{R} \times (-\frac{d}{2}, \frac{d}{2})$
are available in $\mathbb{R}^n \times (-\frac{d}{2}, \frac{d}{2})$ for $n \ge 2$ if we
adapt the regularity properties of the initial and boundary conditions.
\end{remark}
\end{document}